\documentclass[10pt]{amsart}
\usepackage{a4}
\newtheorem{thm}{Theorem}[section]

\newtheorem{lem}[thm]{Lemma}
\newtheorem{prop}[thm]{Proposition}
\theoremstyle{definition}
\theoremstyle{definition}
\newtheorem{ex}{Example}[section]\theoremstyle{definition}

\theoremstyle{remark}
\newtheorem{rem}[thm]{Remark}
\numberwithin{equation}{section}

\newcommand{\R}{\mathbb R}
\newcommand{\E}{\mathbb E}
\newcommand{\N}{\mathbb N}
\def\tr{\,\mathrm{tr}\,}
\def\P{\mathbb P}
\def\1{\mathbb I}
\def\I{\mathcal I}
\def\J{\mathcal J}
\def\L{\mathcal L}
\def\V{\mathcal V}

\def\b{\beta}

\def\e{\epsilon}

\def\l{\lambda}

\def\n{\nu} 

\def\s{\sigma}

\def\t{\tau}
\def\eps{\varepsilon}

\def\half{\frac{1}{2}}
\def\ov{\overline}

\def\dz{\Delta z}

\def\bs{\tilde\sigma}
\def\bb{\tilde b}
\def\bc{\tilde c}
\def\bff{\tilde f}
\def\bbeta{\tilde \eta}
\def\bu{\tilde u}
\def\bx{\bar x}
\def\by{\bar y}

\def\alp{\alpha}
\def\ra{\rightarrow}

\begin{document}
\title[Finite element  scheme for integro-partial differential equations]{A finite element like scheme for integro-partial differential
  Hamilton-Jacobi-Bellman equations}
\author {Fabio Camilli and Espen R. Jakobsen}
\address{Dip. di Matematica Pura e Applicata, Univ. dell'Aquila, loc. Monteluco di Roio, 67040  l'Aquila,
ITALY} \email{camilli@ing.univaq.it}
\address{Department of Mathematical Sciences,
Norwegian University of Science and Technology, 7491 Trondheim,
NORWAY} \email{erj@math.ntnu.no}
\date{\today}

\begin{abstract}
We construct a finite element like scheme for fully non-linear
integro-partial differential equations arising in optimal control of
jump-processes. Special cases of these equations include optimal
portfolio and option pricing equations in Finance. The schemes are
monotone and robust. We prove that they converge in very general
situations, including degenerate equations, multiple dimensions,
relatively low regularity of the data, and for most (if not all)
types of jump-models used in Finance. In all cases we provide
(probably optimal) error bounds. These bounds apply when
grids are unstructured and integral terms are very singular, two
features that are new or highly unusual in this setting.
\end{abstract}
\thanks{The research was supported by the Research Council of Norway
  through the project ``Integro-PDEs: Numerical methods, Analysis, and
  Applications to Finance".}

 \maketitle
\section{Introduction}
In this paper we introduce and analyze finite element (FEM) like
schemes for nonlocal Hamilton-Jacobi-Bellman equations (HJB
equations) of the form
\begin{align}
    \sup_{v\in V}\left\{-\tr [a(x,v)D^2 u]- b(x,v)Du+
    c(x,v)u-f(x,v)-\I^v u(x)\right\}=0\quad\text{in}\quad\R^N,\label{HJBF}\\
\intertext{or}
\label{HJB}
    \sup_{v\in V}\left\{-\tr [a(x,v)D^2 u]- b(x,v)Du+
    c(x,v)u-f(x,v)-\J^v u(x)\right\}=0\quad\text{in}\quad\R^N,
\end{align}
where
\begin{align}
  a(x,v) &=\half \sigma(x,v)\sigma^T(x,v),\\
  \I^v u(x) &=\int_E[u(x+\eta(x,v,z))-u(x)]\nu(dz), \label{IDop}\\
\J^v u(x) &=
  \int_E[u(x+\eta(x,v,z))-u(x)-1_{|z|<1}\eta(x,v,z)Du(x)]\nu(dz),\label{JDop}
\end{align}
$V$ is a compact metric space  and $E=\R^M\setminus\{0\}$. The
coefficients $\sigma,\eta,b,c,f$ are loosely speaking Lipschitz
continuous in $x$, continuous in $v$, and Borel measurable in $z$.
The precise assumptions will be given in the next section. These
equations are the dynamic programming equations
for stochastic control problems involving  Levy processes, a class
of Markov processes with jumps \cite{BJK:08}. The two forms of the equation
correspond  to different intensities $\nu(dz)$ of small jumps in the
corresponding Levy processes. In general these equations are
degenerate and fully-nonlinear, and the corresponding solutions are
typically only H\"older continuous and have to be
understood in the viscosity sense \cite{FS:Book,JK03,BI07}.

Equations like \eqref{HJBF} and \eqref{HJB} appear in advanced
models of financial markets where price evolution of stocks (and
other risky assets) are modeled as (exponential) pure-jump or
jump-diffusion processes. Special cases of \eqref{HJBF} and
\eqref{HJB} include linear equations used in option pricing
problems, obstacle problems used e.g. in pricing of American
options, and fully non-linear equations (the full HJB equation) used
in optimal portfolio problems. We refer to \cite{ContTankov:Book}
for the first two cases and to \cite{FS:Book,BKR} for the last case.
It is well-known that the standard Black-Scholes model (a diffusion
model) give poor fit to real markets, at least on smaller time
scales. E.g. log-returns distributions of stock prices are
leptokurtic and have longer and fatter tails than predicted by the
Black-Scholes model, see e.g. \cite{ContTankov:Book}. To improve
upon these shortcomings, many pure-jump and jump-diffusion models
have been suggested in the literature over the years, see e.g.
\cite{ContTankov:Book,Ry97} for the most popular models. The
empirical fact that Levy processes with discontinuous sample paths
tend to better model e.g. stock prices, is one main reason for the
renewed interest in stochastic optimal control of jump-diffusion
processes.

Except in very simple cases equations \eqref{HJBF} and
\eqref{HJB} do not have closed form solutions, so numerical methods
are needed to obtain (approximate) solutions. In this paper we
construct FEM like schemes. Since the equations are fully nonlinear
we have no weak/variational formulation, and the usual FEM approach
does not work. 
If the measure $\nu$ is finite and the corresponding HJB
equation is \eqref{HJBF}, we discretize in two steps. First we obtain
a semi-discrete equation (``discrete in time'') as the dynamic
programming equation  of a discrete time control problem approximating
the underlying continuous time control problems of \eqref{HJBF}:
\begin{align}
u_h(x)&=\inf_{v\in V }\Big\{hf(x,v)+ e^{-h c(x,v)}\Big[\frac{e^{-h \l
  }}{2d}\sum_{m=1}^d \Big(u_h(x+h b(x,v)+\sqrt{h}\s_m(x,v))+ \label{S1}
\\
&+u_h(x+h b(x,v)-\sqrt{h}\s_m(x,v)\Big)+
\frac{1-e^{-h \l}}{\l}\int_E u_h(x+\eta(x,v,z))\nu(dz)\Big]\Big\},\nonumber
\end{align}
where $\sigma_m$ is the $m$-th column of the matrix $\sigma$, $h$ is
the discretization step and $\l$ is the mass of the measure. The
next step is obtain a fully discrete equation by introducing a regular
triangulation and look for continuous piecewise linear functions over
the chosen triangulation satisfying the semi-discrete equation at
every vertex of the triangulation.
When $\nu$ is not finite we first approximate it by a finite measure
$1_{r<|z|<R}\nu(dz)$ (truncation), and then approximate the truncated
equation following the above approach. To improve the truncation
approximation, we also add small diffusion and/or drift terms to the equation.

We prove that these methods converge and derive
(probably optimal) error bounds using the framework developed in
\cite{JKLC} and ideas from \cite{CF}. We also discuss issues like
restricting to a bounded domain, truncating long jumps, and
approximating integrals by quadrature. In all cases we provide
rigorous error bounds for the various approximations. What remains to do in a
computer implementation, is the resolution of the non-linearity. There
are various ways to do that, e.g. policy iteration, value iteration,
artificial time methods (stationary limits), and so on. We will not
address this point in this paper.

FEM  schemes like the one we have described above are usually called
semi-Lagrangian or control schemes in the literature. They are
(usually) monotone and first order accurate, and a comprehensive
background and references can be found in \cite{F}.  Most results in
this field concerns deterministic control problems and first order
HJB equations without integral terms. Semi-Lagrangian schemes for
second order HJB equations with no integral terms  have been considered in
\cite{Me:Estim,CF,BJ:Err1,AO}. Moreover in \cite{Ca:IPDE} such schemes
were derived for HJB equations associated to piecewise deterministic
processes (compound Poisson processes with drift). One advantage of
this type of schemes is that in general they produce also an
approximation of the optimal control law in feedback form and the
optimal trajectories (this is a key point in the study of a control
problem). This advantage can also be shared by some (monotone, low
order) finite difference schemes following the construction of
Kushner, and we refer to \cite{KuDu} for a discussion of this point.
An other advantage is that these schemes are formulated on general
grids/triangulations just like FEMs. By contrast this is very
cumbersome to achieve with finite difference methods.  A third
advantage is that the semi-Lagrangian scheme {\em automatically}
handles non-diagonally dominant diffusion matrices $a$. Such
matrices appear in applications and are cumbersome to discretize, we
refer to \cite{BZ} for a finite difference (FDM) approach to this
problem.

The construction and analysis of numerical schemes for linear
integro-partial differential equations arising as pricing equations
in financial markets of jump-diffusion type is currently an active
field of research, see e.g.
\cite{ContTankov:Book,MvPS:04,dHFL:04,BNR07} and references therein.
By contrast, there are few works on numerical schemes for fully
nonlinear degenerate integro-partial differential equations. We
mention the discussion about jump-diffusion processes in \cite{KuDu}
and the papers \cite{Said-KHK,BLCN,JKLC,BJK2}. In all cases mainly
monotone FDMs are considered and convergence is obtained. In the
last two papers the main focus is on convergence rates. The
framework of \cite{JKLC} is again based on ideas of Krylov, Barles
and Jakobsen \cite{Kry97,Kry00,BJ:Err1} in the pure PDE case. We
emphasize that the error bounds of this paper apply when grids are
unstructured and integral terms are very singular, a feature which
is new or highly unusual in this setting.

The rest of this paper is organized as follows. In section
\ref{sec:prelim} we state the assumptions on the data and give
well-posedness/regularity results for equations \eqref{HJBF} and
\eqref{HJB}. We discuss truncation of the Levy measure and reduction
to a bounded domain and give error bounds. In section
\ref{sec:constr} we construct the schemes via the dynamic
programming equation (a semi-discretization) and FEM ideas. We prove
existence, uniqueness, consistency, and partial error bounds. We
also derive a fully discrete scheme of FEM type by piecewise linear
reconstruction on a (regular) triangulation of the domain. In
section \ref{sec:conv} we derive error bounds for the semi-discrete
scheme with or without truncation of the Levy measure. Finally in
the appendix, we prove the main technical result of this paper, a
regularity and continuous dependence result for the semi-discrete
scheme.





\subsection*{Notation:}
By $USC_b(\R^N)$ and $C_b(\R^N)$ we mean the spaces of bounded
continuous and upper semicontinuous functions.
We will use the following norms
\begin{eqnarray*}
    |f|_0=\mathrm{ess \,sup}_{x\in\R^N}|f(x)|, \quad
    [f]_1=|D f|_0, \quad \mbox{and} \quad
    |f|_1=|f|_0+[f]_1.
\end{eqnarray*}

\section{Preliminaries}
\label{sec:prelim}

In this section we state the assumptions for equations \eqref{HJBF}
and \eqref{HJB} and give well-posedness results. We discuss reduction
to a bounded domain and reduction to bounded measure $\nu$ with compact
support.

\subsection{Assumptions, well-posedness, regularity}
We will use the following assumptions:

\begin{itemize}
\item[(A1)] The set $V$ is a compact metric space, the coefficients
  $\sigma, \eta, b, c, f$ are continuous in $x$ and $v$ and Borel
  measurable in $z$, and $\nu$ is a positive Radon measure on $E:=\R^M\setminus\{0\}.$\\

\item [(A2)] There exists $L_1,L_2, \ell \geq0$ such that for any
  $v\in V$ and $z\in E$
    \begin{align*}
      &|c(\cdot,v)|_1 +|f(\cdot,v)|_1\leq L_1,\\
      &|\sigma(\cdot,v)|_1
      +|b(\cdot,v)|_1 \leq L_2\\
      &|\eta(\cdot,v,z)|_1\leq L_2( |z|1_{|z|<1}+ e^{\ell|z|}1_{|z|>1}).
    \end{align*}

\item [(A3)] There exists $c_0>0$ such that, for
all $x\in\R^N$ and $v\in V$,
$$
  c(x,v)\geq c_0.
$$

\item[(A4)] There exists $L_3\geq0$, such that for $\ell$ defined in
  (A2), any $v\in V$ and $x\in \R^N$,
 $$|D_z\eta(x,v,\cdot)|_0\leq L_3 e^{\ell|z|}.$$
\end{itemize}
\medskip

For the measure $\nu$ we will use the following integrability
assumptions.
\begin{itemize}
\item[(B1)$_k$] The measure $\nu$ satisfies for $\ell$ defined in (A2)
  and some $k\in\{0,1,2\}$,
    $$
      \int_{E}(|z|^k 1_{|z|<1} +e^{\ell |z|}1_{|z|>1})\,\nu(dz)<\infty.
    $$



 \item [(B2)] The measure $\nu$ has a positive density $m:E\to
    [0,\infty)$ such that $\nu(dz)=m(z)dz$ and for some $C$, $\eps>0$,
    and $n\in\N$
    \begin{align*}
       |D^km(z)|\leq C_ke^{-(\ell+\epsilon)|z|}\quad\text{for}\quad
     k=0,1,\dots,n ,
    \end{align*}
 where $D^km$ is the vector of all order $k$ derivatives of $m$ and
$\ell$ is defined in (A2).
\medskip

 \item [(B3)] The measure $\nu$ has a positive density $m:E\to
    [0,\infty)$ such that $\nu(dz)=m(z)dz$ and for some
    $\alp\in[0,2)$, $C_0$, $C_1$, and $\eps>0$,
    \begin{align*}
       &|D^km(z)|\leq C_k\Big( \frac
     1{|z|^{N+\alpha+k}}
     1_{|z|<1}+e^{-(\ell+\epsilon)|z|}1_{|z|>1}\Big)\quad\text{for}\quad
     k=0,1 ,
\end{align*}
where $D^0m=m$, $D^1m=Dm$ the gradient of $m$,
and $\ell$ is defined in (A2).
\end{itemize}


Assumptions (A1) -- (A3) and (B1)$_2$ are standard from a stochastic control
theory point of view, as they insure existence and uniqueness of
strong solutions of the underlying stochastic differential equations,
see \cite{GiSk}.
Under assumption (B1)$_2$ or the less general assumption (B1)$_1$ the
measure $\nu$ may have a (non-integrable) singularity at
$z=0$. 
If (B1)$_1$ or (B1)$_0$ holds the HJB equation takes the form \eqref{HJBF}. If
(B1)$_2$ holds but not (B1)$_1$, then \eqref{HJB} gives the correct
form of the HJB equation. In this case, the extra term in the
integrand of $\J^v$ (compared to $\I^v$) is needed for the integral to
converge: If (B1)$_2$ holds, then $\J^v u$ converge for all
$C^2$ functions $u$ is with polynomial growth at infinity. Finally,
(B2) and (B3)  prescribe densities.
 Close to $z=0$, the density of (B3) equals the density of the
 $\alp$-stable processes related to the fractional Laplacian
 $\Delta^{\alp/2}$. (B3) implies (B1)$_2$, and if
$\alp\in[0,1)$, then (B3) also implies (B1)$_1$.

\begin{ex}
When the stock returns are modeled as exponential Levy processes the
 integral terms equals \eqref{IDop} or \eqref{JDop} with
 $\eta(x,z)=x(e^z-1)$ in the  one-dimensional uncontrolled case, see
 \cite{ContTankov:Book}. An
 important extension of this expression satisfying (A2) and (A4) is
$\eta(x,v,z)=\bar\eta(x,v)\phi(z)$ where $\bar \eta$ is bounded and
 Lipschitz and $|\phi(z)|\leq C(e^{|z|}-1)$, $|D\phi(z)|\leq Ce^{|z|}$.
\end{ex}

\begin{ex}
Assumptions (B2) and (B3) is satisfied by all bounded and unbounded
L\'evy processes used in the literature to model financial markets,
see \cite{ContTankov:Book} for a nice overview. In the models of
Merton and Kou,
the Levy measures have bounded densities (no singularity), in the
Merton case \cite{Merton} given by
$$\nu(dz)=\frac{\lambda}{\delta\sqrt{2\pi}}
e^{-\frac{|z-\mu|^2}{2\delta^2}}dz \quad \text{for constants
  $\lambda,\delta,\mu$}.$$
For these models (B2) holds. The Variance Gamma model has a Levy
density with an integrable singularity at $z=0$, and satisfy (B3)
with $\alp=0$. The Normal Inverse Gaussian model \cite{Ry97} has a
non-integrable density corresponding to $\alp=1$ in (B3), and the
Levy measure is
$$\nu(dz)=\frac{C}{|z|}e^{Az}K_1(B|z|) \quad \text{for constants
  $A,B,C$},$$
where $K_1$ is a modified Bessel function of 2nd kind. Finally we
  mention models using tempered $\alp$-stable processes, e.g. the
  CGMY model \cite{CGMY}. Here $\alp\in(0,2)$ in
  (B3) and the Levy measure is give by
\begin{align*}
\nu(dz)=\frac{c_-}{|z|^{1+\alp}}e^{-\lambda_-|z|}1_{z<0}+
  \frac{c_+}{|z|^{1+\alp}}e^{-\lambda_+|z|-}1_{z>0}\\ \text{for
  constants
  $\alp\in(0,2),c_-,c_+,\lambda_-,\lambda_+\geq0$.}
\end{align*}
\end{ex}

It is well known that under the above assumptions the solutions to
\eqref{HJBF} and \eqref{HJB} need not be smooth, and that the correct
concept of (weak) solutions is that of viscosity
solutions. For the definition of viscosity solution in this case we refer to
\cite{Sa:IPDE1,BBP97,Pham,JK03,BI07}. We state without proof a
well-posedness and regularity result for \eqref{HJB}. The proof of
this result is standard, and we refer to
\cite{BBP97,Pham,JK:CDIPDE,JK03} for the proofs of similar results.

\begin{thm}
\label{WP}
Assume (A1) -- (A3) and (B1)$_2$ hold.

(i) There exists a unique viscosity solution $u\in C_b(\R^N)$ of equation
\eqref{HJB} which is H{\"o}lder continuous, i.e., there is
a $\delta\in(0,1]$ such that
$$
|u(x)-u(y)|\leq C|x-y|^\delta \quad \text{for all}\quad x,y\in\R^N.
$$

(ii) There exists a constant $c_1>0$ depending only on
$\sup_v[\sigma(\cdot,v)]_1$, $\sup_v[b(\cdot,v)]_1$, and \linebreak
$\sup_v\int_E[\eta(\cdot,v,z)]_1^2\nu(dz)$
such that if $c_0 \geq c_1$, then the viscosity
solution $u$ of \eqref{HJB} is Lipschitz continuous ($\delta=1$ above).
\medskip

(iii) Let $u,-v\in  USC_b(\R^N)$. If $u$ and $v$ are
respectively viscosity sub- and supersolutions of \eqref{HJB}, then
$u\leq v$ in $\R^N$.
\end{thm}

\begin{rem}
\label{remWP} Note that this result also holds for \eqref{HJBF}
under assumptions (A1) -- (A3) and either one of (B1)$_1$ or
(B1)$_0$. Simply write this equation in the form \eqref{HJB} with
$b$ replaced by $\bar b= b-\int_E \eta\nu(dz)$.
\end{rem}

\subsection{Reduction to bounded measure with compact support}

Consider Levy measures $\nu$ which are not bounded nor have compact support.
We assume that (B3) hold and consider two cases: (i) $\alp\in[0,1)$
  and (ii) $\alp\in[1,2)$. In the first case we truncate the Levy and
  compensate by adding a drift and a diffusion term. Let us introduce
  the ``two-scales'' truncated L\'evy measure
\begin{equation}\label{trunc}
\nu_{r,R}(dz):=1_{r<|z|<R}\,\nu(dz),
\end{equation}
where $r\in (0,1)$ and $R>1$. Clearly, $\nu_{r,R}(dz)$ is a bounded
and compactly supported measure satisfying assumption (B1)$_0$. When
$\alp\in(0,1)$ assumption (B1)$_1$ holds and the HJB equation is
\eqref{HJBF}. This equation is approximated by
\begin{align}
    \sup_{v\in V}\left\{-\tr [\bar a(x,v)D^2 u]- \bar b(x,v)Du+
    c(x,v)u-f(x,v)-\bar \I^v
    u(x)\right\}=0\quad\text{in}\quad\R^N,\label{T1HJB}
\end{align}
where
\begin{align}
\bar a(x,v) &=a(x,v)+\frac12\int_{0<|z|< r}
\eta(x,v,z)\eta(x,v,z)^T \,\nu(dz),\label{bara}\\
\bar b(x,v) &=  b(x,v)+\int_{0<|z|<r}\eta(x,v,z)\nu(dz),\label{barb}\\
\bar \I^v
    u(x) &=\int_{E}\left[\phi(x+\eta(x,v,z))-\phi(x)
    \right]
    \,\nu_{r,R}(dz).\label{apprI}
\end{align}
This gives a ``third order'' approximation of $\I^v\phi$ since
$$|\bar
\I^v\phi+(\bar b -b)D\phi+\tr[(\bar a-a)D^2\phi] - \I^v\phi|\leq
K|D^3\phi|_0\left|\int_{0<|z|<r}|\eta(\cdot,\cdot,z)|^3\nu(dz)\right|_0.$$

In case (ii) we truncate the Levy measure and compensate by
adding a diffusion term. When $\alp\in[1,2)$ assumption (B1)$_1$ is not
  satisfied and the HJB equation is \eqref{HJB}. We approximate this
  equation by
\begin{align}
    \sup_{v\in V}\left\{-\tr [\bar a(x,v)D^2 u]- \tilde b(x,v)Du+
    c(x,v)u-f(x,v)-\bar \I^v
    u(x)\right\}=0\quad\text{in}\quad\R^N,\label{THJB}
\end{align}
where $\bar a$ and $\bar \I^v u$ is defined in \eqref{bara} and
\eqref{apprI}, and
\begin{align}
\label{tildeb}
\tilde b(x,v) &=  b(x,v)+\int_E 1_{|z|<1}\eta(x,v,z)\nu_{r,R}(dz).
\end{align}
Again we obtain a ``third order'' approximation of the integral term
($\J^v\phi$ this time) since
$$|\bar
\I^v\phi+(\tilde b -b)D\phi+\tr[(\bar a-a)D^2\phi] - \J^v\phi|\leq
K|D^3\phi|_0\left|\int_{0<|z|<r}|\eta(\cdot,\cdot,z)|^3\nu(dz)\right|_0.$$

Equations \eqref{T1HJB} and \eqref{THJB} are almost, but not
quite, of the same form as \eqref{HJBF}. In \cite{JKLC} it is proved
that Theorem \ref{WP} still holds for solutions of \eqref{T1HJB} and
\eqref{THJB}. The following result gives error bounds for the
approximations \eqref{T1HJB}, \eqref{THJB} of equations \eqref{HJBF},
\eqref{HJB}.

\begin{lem}
\label{TruncThm} Assume (A1) -- (A3) and (B3) hold and let $u, \bar u,
v, \bar v \in C^{0,1}(\R^N)$ solve \eqref{HJBF}, \eqref{T1HJB},
\eqref{HJB}, \eqref{THJB} respectively. Then for $r>0$ small enough and $R$
large enough we have
$$|u- \bar u|_0+|v- \bar v|_0\leq C_1(r^{1-\alp/3}+e^{-\ell R}),$$
for some constant $C_1$ independent of $r$ and $R$.
\end{lem}

The result about $v,\bar v$ is proved in \cite{JKLC}, while the result
about $u,\bar u$ is new but the proof is almost identical to the proof
 of the $v,\bar v$ result. We omit it.

\begin{rem}
The compensating drift and diffusion terms are added to improve the
convergence of the approximation (of the small jumps). The error is
of order $r^{1-\alp/2}$ without the the compensating diffusion (see
\cite{JKLC}), while it is of order $r^{1-\alp}$ without compensating
drift and diffusion in the case $\alp\in[0,1)$.
\end{rem}

\begin{rem}
In probabilistic terms, the explanation for the improved convergence
is that the small jumps of a Levy process can be approximated by a
Brownian motion, we refer to e.g. \cite{AsRo,ContTankov:Book} for
details. Moreover using probabilistic methods (Berry-Essen type
estimates), it is possible to prove convergence of order $r$ in some
cases \cite{AsRo}. We refer to \cite{CV:FDM} where such estimates
are made explicit for linear parabolic equations with constant
coefficients.
\end{rem}

\begin{rem}
When $\alp=0$  no improvement is obtained by  adding drift and
diffusion terms to the equation. If $\nu$ has a bounded
density, no truncation is needed  close to $z=0$.
\end{rem}

\subsection{Reduction to a bounded domain}
\label{sec:BndDom} Reduction to a bounded domain is key step in
order to implement a numerical method. Following ideas of \cite{CF},
we  restrict equations \eqref{HJBF} and \eqref{HJB} to bounded
domains by truncating the coefficients outside some large ball,
\[
    B_{\frac1\mu} =\{ x\in  \R^n  :|x|<{\frac{1}{\mu}}\}\quad\text{for}\quad \mu>0.
\]
Let $\xi_\mu\in C_c^\infty(\R^N)$ be a cut-off function satisfying
$0\le \xi_\mu\le 1$ and $\xi_\mu(x) = 1$ for $x \in
B_{\frac1\mu}$ and define
\[
    \sigma_\mu= \xi_\mu (x)\sigma(x,a),\qquad b_\mu =\xi_\mu (x)
    b(x,a),\qquad \eta_\mu(x,v,z) =\xi_\mu (x) \eta(x,v,z).
\]
Since $b_\mu$, $\sigma_\mu$ and $\eta_\mu$ satisfies the same
assumptions of $\sigma$, $b$ and $\eta$ for any $\mu$, there exists a
unique viscosity solution $u_\mu$ of the equation
\begin{equation}\label{HJBm}
    \sup_{v\in V}\left\{-\tr [a_\mu(x,v)D^2 u]- b_\mu(x,v)Du+
    c(x,v) u-f(x,v)-\J^v_\mu(x)\right\}=0\quad\text{in}\quad\R^N\\
\end{equation}
where $ a_\mu(x,v) =\half \sigma_\mu(x,v)\sigma_\mu^T(x,v)$ and
$\J^v_\mu (\cdot)$ is defined as in \eqref{JDop} with $\eta_\mu$ in
place of $\eta$.  Since the coefficients $b_\mu$, $\sigma_\mu$ and
$\eta_\mu$ are zero outside of ${\rm supp}(\xi_\mu)$, the solution $u_\mu$
\eqref{HJBm} is given by
\[u_\mu(x)=\min_{v\in V}\frac{f(x,v)}{c(x,v)} \quad \text{for
  any}\quad x\not\in \mathrm{supp}(\xi_\mu).\]

The next result give a crude bound on the error that this cut-off
procedure introduces.
\begin{lem}
    Assume (A1), (A2), (B1)$_2$ hold, that \eqref{HJB} and \eqref{HJBm}
  satisfy the dynamic programming principle, and that $u$, $u_\mu$
  solve \eqref{HJB}, \eqref{HJBm} respectively. Then there exists a
  constant $C$ such that
\[|u(x)-u_\mu (x)|_0\leq
    C \mu^2(1+|x|^2)  \quad\text{in}\quad B_{\frac1\mu}.  \]
\end{lem}
We skip the proof since it is similar to the proof of Proposition
3.1 in \cite{CF} when we are equipped with the moment estimates
(3.2) and (3.4) of \cite{Pham}. The dynamic programming principle
has recently been extended to the current setting in \cite{BJK:08}.

\begin{rem}
As in Remark \ref{remWP} we immediately get an analogous result for
\eqref{HJBF} under assumption (B1)$_1$.
\end{rem}

\begin{rem}
If the equation is uniformly elliptic or has a non-degenerate
singular integral term, then you expect estimates
decaying exponentially as $\mu\ra\infty$. We refer to
\cite{CV:FDM,MvPS:04} for such result in a linear one-dimensional
setting.
\end{rem}

\section{Construction of the scheme}
\label{sec:constr}

In subsections \ref{sec:timediscr} and \ref{sec:FullyDiscr} we  always
assume that (A1) -- (A3)
and (B1)$_0$ hold: The Levy measure $\nu$ is bounded and the HJB
equation has the form \eqref{HJBF}. We can always reduce to this case
by a truncation. In subsection
\ref{sec:unbnd} we discuss the general case.

\subsection{Semi-discretization}
\label{sec:timediscr}

We introduce a control problem for which \eqref{HJBF} is the
corresponding dynamic programming or HJB equation. We start by
defining the controlled dynamics. Since the measure $\nu$ is
bounded, we can normalize $\nu$ and obtain a probability measure
$\mu$ as follows
\[\mu(dz)=\frac{1}{\l}\nu(dz)\]
where
\begin{equation}\label{rate}
  \l=\int_E\nu(dz)
\end{equation}
Now we consider a Markov process $X_t^{v_\cdot}$ evolving  according  the SDE
\begin{equation}\label{SDE}
    dX_t= b(X_{t^-},v_{t^-})dt +
    \sigma(X_{t^-},v_{t^-})dW_t+\int_{|z|>0}\eta(X_{t^-},v_{t^-},z)
    \bar{\mu}(dz,dt),
\end{equation}
where $\bar{\mu}$ is a Poisson measure corresponding to a compound
Poisson process with jump intensity $\l$ and jump distribution
$\mu$. The control $v_\cdot$  belongs to $\V$, the set of all
progressively measurable processes with values in $V$. Since we
assume that $\nu$ and hence $\bar\mu$ are bounded, on any finite
time interval $X_t$ will only jump finitely many times with
probability one. Between two jump times $T_i$ and $T_{i+1}$ the
process diffuses according to the SDE
\begin{equation}\label{diffSDE}
    dX_t= b(X_{t^-},v_{t^-})dt + \sigma(X_{t^-},v_{t^-})dW_t.
\end{equation}
For $v_t\equiv v\in V$ and $a(x,v)=\frac12\sigma(x,v)\sigma^T(x,v)$, the
infinitesimal generator of the process $X_t$ is
\begin{align*}\label{IG}
    L^v\psi(x)&=\lim_{t\to
    0}\frac{\E_x[\psi(X_t)]-\psi(x)}{t}=\mathrm{tr}\, [a(x,v)D^2
    \psi]+ b(x,v)D\psi+\I^v \psi(x)
\end{align*}
for $\psi\in C^2(\R^N)$. On the paths of the  process $X_t$ we
define the discounted cost functional
\[
    J(x,v_\cdot)=\E_x\Big[\int_0^\infty f(X_t,v_t)e^{-\int_0^t
    c(X_s,v_s)ds}dt\Big], \]
and we consider the corresponding value function
\[
u(x)=\inf_{v_\cdot\in\V} J(x,v_\cdot).
\]
In \cite{BJK:08} it is proved that $u$ is the unique viscosity solution
of equation \eqref{HJBF}.

Following the approach of \cite{Ca:IPDE,CF} we construct an
approximation scheme for the equation \eqref{HJBF} by discretizing
the associated control problem. We fix a discretization step $h>0$
and consider two stochastic processes $N_n$ and $Z_n$, $n\in \N$,
taking values in $\N$ and in $\R^N$ and representing the $n$-th jump time and
the corresponding $z$-jump (size and direction) of the Poisson
measure $\bar \mu$. We set $N_0=0$ and $Z_0=0$
and assume that $N_n$ has independent
$h\l$-exponentially distributed increments, i.e. the probability
distribution of $N_n$ is given by
\begin{align*}
\P[N_{n+1}-N_n\ge j| N_0,N_1,\dots,N_n]=e^{-h\l j },\quad n=0,1,2,\dots,
\end{align*}
while the $Z_n$, $n\in\N$, are i.i.d. random variables with
probability density $\mu$. Now we define a discrete time stochastic
process $X_n$ approximating the continuous time process
$X_t$.
\begin{align*}
\begin{cases}
X_0=x,\\
X_{n}=X_{n-1}+h b(X_{n-1},v_{n-1})+\sqrt{h} \sum_{m=1}^d\s_m(X_{n-1},v_{n-1})\xi^m_{n-1},\hspace{3cm}\\
\qquad$\hfill$\text{for $n=N_i+1,N_i+2,\dots, N_{i+1}-1$,}\\
X_{N_{i+1}}=X_{N_{i+1}-1}+\eta(X_{N_{i+1}-1},v_{N_{i+1}-1},Z_i),
\end{cases}
\end{align*}
for $i=1,2,3,\dots$, where $\sigma_m$ denote the $m$-th column of
$\sigma$ and $\xi^m_n$, $m=1,\dots,d$ are random variables taking
values in $\{-1,0,1\}$ such that
$$\P[\{\xi_n^i=\pm1\}]=\frac1{2d}\qquad\text{and}\qquad
  \P[\{\xi_n^i\neq 0\}\cap\{\xi_n^j\neq0\}]=0, \quad i\neq j.$$
The discrete control $\{v_n\}$ is a random variable with values in $V$
which is measurable with respect to the $\s$-algebra generated by
$X_1,\dots,X_n$.

Between jumps the process evolves like a random walk approximating
the SDE \eqref{diffSDE}, and when the process jumps there is no
diffusion/random walk. The  generator of the discrete process is
\begin{equation}\label{IG}
    L_h^v\psi(x)=\frac{\E_x[\psi(X_1)]-\psi(x)}{h}= e^{-h\l}\L_h^v \psi(x)
               +\frac{1-e^{-h\l}}{h}\I^v_h \psi(x),
\end{equation}
for $\psi\in C^0(\R^N)$ and where
\begin{align*}
     &\L_h^v\psi=\frac{1}{2dh}\sum_{m=1}^d\Big[\psi(x+h b(x,v)+\sqrt{h}\s_m(x,v))+\psi(x+h b(x,v)-\sqrt{h}\s_m(x,v))-2\psi(x)\Big],\\
      &\I_h^v\psi=\frac{1}{\l}\I^v\psi,
\end{align*}
with $\I^v$ as in \eqref{IDop}. Observe that at this level the space
variable is not discretized, therefore the discrete process has the
same jump distribution as the continuous process.

On the paths of the discrete process we define the cost functional
\begin{equation}\label{cost}
    J_h(x,\{v_n\})=\E_x[\sum_{n=0}^\infty he^{-h\sum_{i=0}^{n-1}
    c(X_i,v_i)}f(X_n,v_n)],
\end{equation}
(with the convention $\sum_{i=0}^{-1}=0$), and the corresponding value
function for the discrete control problem
\begin{align}
\label{VF}
u_h(x)=\inf_{\{v_n\} } J_h(x,\{v_n\}).
\end{align}
Now it is easy to see, at least formally  \cite{BSh}, that  the following
dynamic programming principle holds:
\begin{equation*}
    u_h(x)=\inf_{\{v_i\} }\E_x\Big[\sum_{n=0}^P he^{-h\sum_{i=0}^{n-1}
    c(X_i,v_i)}f(X_n,v_n)+e^{-h\sum_{i=0}^{P} c(X_i,v_i)}u_h(X_{P+1})\Big],
\end{equation*}
for any $P\in\N$. Taking $P=0$ in the above equation and noting that
$X_0=x$, we get
\begin{align}
&u_h(x)=\inf_{v\in V }\E_x[hf(X_0,v)+ e^{-h c(X_0,v)}u_h(X_1)]\nonumber\\
&=\inf_{v\in V }\Big\{hf(x,v)+
 e^{-h c(x,v)}\Big[\frac{e^{-h \l }}{2d}\sum_{m=1}^d (u_h(x+h
   b(x,v)+\sqrt{h}\s_m(x,v))\label{DPPh}\\
&+u_h(x+h b(x,v)-\sqrt{h}\s_m(x,v))+
\frac{1-e^{-h \l}}{\l}\int_E u_h(x+\eta(x,v,z))\nu(dz)\Big]\Big\}.\nonumber
\end{align}
Rearranging the terms in  the previous equation and dividing by
$he^{-hc(x,v)}$ we get
\begin{equation}\label{HJBFh}
    \sup_{v\in V}\left\{-L_h^v u_h(x)+
    \frac{e^{hc(x,v)}-1}{h} u_h(x)-e^{hc(x,v)}f(x,v)\right\}=0 \quad
    \text{in}\quad\R^N,
\end{equation}
where $L_h^v(\cdot)$ is as in \eqref{IG}. We will talk about sub-
and supersolutions of this equation, meaning that \eqref{HJBFh}
holds as an inequality with $\leq$ and $\geq$ respectively. For the
scheme \eqref{HJBFh} we have the following easy properties:

\begin{prop}\label{propscheme1}
Assume (A1) -- (A3) and (B1)$_0$.
\smallskip

(i) If $u_h$ and $v_h$ are bounded sub- and
supersolutions of \eqref{HJBFh}, then $u_h\leq v_h$ in $\R^N$.
\smallskip

(ii) Any solution $u_h$ of \eqref{HJBFh} is bounded and satisfies for
all $h>0$,
$$ |u_h |_0\le \frac{\sup_v|f(\cdot,v)|_0}{c_0},\quad\text{where $c_0$ is as in (A3)}.$$

(iii) There exists a unique bounded continuous function $u_h$
solving \eqref{HJBFh}.
\smallskip

(iv) The solution of \eqref{HJBFh} is (at least formally) equal to
the value function \eqref{VF}.
\smallskip

(v) For $0<h<1$ and $\phi\in C^4(\R^N)$ satisfying $|\phi|_0+
\dots+|D^4\phi|_0<\infty$,
$$|L^v\psi(x) - L_h^v\psi(x)|\leq
C_1h(|D^2\phi|_0+|D^3\phi|_0+|D^4\phi|_0)+C_2h\lambda((1+|\int_E\eta\nu|)|D\phi|_0+|D^2\phi|_0),$$
where the constants $C_1$ and $C_2$ only depend on
$\sup_v|\sigma|_0,\sup_v|b|_0$.
\end{prop}

\begin{rem}
For the truncation of an unbounded measure to converge as $h\ra0$,
we need $\lambda \ra\infty$ as $h\ra 0$, while the
scheme \eqref{HJBFh} converges only if both $h\ra0$ and $\lambda
h\ra 0$ by {\it (v)}. The last condition means that the small jumps must be
resolved in the grid.
\end{rem}

\begin{proof}
We work with the scheme in the equivalent form \eqref{DPPh}. Note that the
scheme is monotone, it has positive coefficients. If $u$ and $v$ be sub-
and supersolutions of \eqref{HJBFh}, an easy and standard computation using
\eqref{DPPh} and assumption (A3) shows that
$$u_h(x)-v_h(x)\leq e^{-c_0h}|(u_h-v_h)^+|_0,$$
and hence $(1- e^{-c_0h})|(u_h-v_h)^+|_0\leq 0$ which proves (i). A
similar computation shows (ii) after noting that $\frac{hc
  e^{hc}}{e^{hc}-1}\leq 1$ for $hc\geq 0$.
Next denote the right hand side of \eqref{DPPh} by $T_h u_h$, and note
that $T_h$ is contraction in the $|\cdot|_0$ norm,
$$T_hu_h-T_hv_h \leq e^{-hc_0}|u_h-v_h|_0.$$
Existence and uniqueness of a continuous bounded solution to
\eqref{HJBFh} follows from Banach's fixed point theorem and this proves
 (iii). Part (iv) follows from the from the dynamic
programming principle as explained before the proposition. Finally, to
prove (v), note that since $\I^v_h=\lambda^{-1}\I^v$ we may write
$$L^v\psi -
L_h^v\psi=(\L^v-\L^v_h)\psi+(\I^v-\I^v)\psi-(1-e^{-h\lambda})\L^v_h\psi-(1-h\lambda^{-1}(1-e^{-h\lambda}))\I^v\psi,$$
where $\L^v\psi(x)=\tr[a(x,v)D^2\psi(x)]+b(x,v)D\psi(x)$.
By Taylor expansion, e.g.
$$|(\L^v-\L^v_h)\psi|\leq
h|b|_0^2|D^2\psi|_0+(h^2|b|_0^3+h
|b|_0|\sigma|_0^2)|D^3\psi|_0+(h^3|b|_0^4+h^2|b|_0^2|\sigma|_0^2
+h|\sigma|_0^4)|D^4\psi|_0,$$ and the estimates
$|1-e^{-x}|,|1-x^{-1}(1-e^{-x})|\leq |x|$ the result follows.
\end{proof}

Next we give an optimal Lipschitz regularity and continuous dependence on
coefficients result for the scheme \eqref{HJBFh}.

\begin{prop}\label{ContLipSch}
Let $u_h$ and $\tilde u_h$ be solutions of \eqref{HJBFh}
corresponding to the data $\sigma, b, c,  f, \eta, \nu$ and $\bs,
\bb, \bc, \bff, \bbeta, \nu$ respectively, assume both sets of
coefficients satisfy (A1) -- (A3) and (B1)$_0$, and that
$h\lambda\leq \bar C_0$ and $h\in (0,1]$. Then there exist constants
$c_1, L, K\geq0$ (only depending on the data and $\bar C_0$) such that if
$c_0\geq c_1$ ($c_0$ as in (A3)), then for all $h>0,x,y\in\R^N$,
\begin{align*}
|u_h(x)-\bu_h(y)|\leq & \ L|x-y| +K\sup_{v\in
V}\Big[|f-\bff|_0+|c-\tilde
c|_0\\
&+|b-\bb|_0+|\sigma-\bs|_0+
\Big|\int_E|\eta(\cdot,z)-\bbeta(\cdot,z)|^2\nu(dz)\Big|_0^{1/2}\Big].
\end{align*}
\end{prop}

If the coefficients are equal and $u_h=\bu_h$, then this is a
Lipschitz regularity result, while if $x=y$ then this is a
continuous dependence on the coefficients result. This result  one of the main
contributions of this paper, and it will play a key role in the next
section where error bounds are derived. It extends similar results of
\cite{BJ:Err1} to equations with integral terms. The proof is rather
long and technical and we have put the main bulk in the appendix. In
the pure diffusion case, the current proof simplifies considerably the
arguments of \cite{BJ:Err1}.

\begin{proof}
The result follows from Theorem \ref{LipCont} in the Appendix after writing
$$b=\Big[b-\frac{1-e^{-h\lambda}}{\lambda e^{-\lambda h}}\int_E\eta
  \nu\Big]+\frac{1-e^{-h\lambda}}{\lambda e^{-\lambda h}}\int_E\eta \nu,$$
and noting that by the Cauchy-Schwartz inequality,
$$\Big|b-\bb - \frac{1-e^{-h\lambda}}{\lambda e^{-\lambda h}}\int_E[\eta-\bbeta]
  \nu\Big|_0\leq
  |b-\bb|_0+
  (e^{h\lambda}-1)\l^{-1/2}\Big|\int_E|\eta-\bbeta|^2\nu\Big|_0^{1/2}.$$
\end{proof}

\begin{rem}
When $c_0<c_1$ the solution to \eqref{HJBFh} is only H\"older
continuous. We will not discuss this case here, but refer instead to
\cite{BJ:Err1} for how to obtain results in this case.
\end{rem}

\subsection{The fully-discrete scheme}
\label{sec:FullyDiscr}
In this section we introduce a FEM like discretization of
\eqref{HJBFh} yielding a fully discrete scheme. For a nice
introduction to FEMs we refer to \cite{BS}.
For $k>0$ let $\mathcal{T}^k=\{S^k_j\}_{j\in\N}$ be a non-degenerate
triangulation of $\R^N$, i.e. a collection of $N$-simplices $S^k_j$ such that
$$\underset{j\in\N}{\cup} S^k_j=\R^N,\quad \sup_{j\in\N}( \hbox{diam}\,S_j)\leq
 k,\quad \rho k\leq \sup_{j\in\N}( \hbox{diam}\,B_{S^k_j}),$$
where  $\rho\in(0,1)$, $\hbox{diam}$ denotes the diameter of the set,
 and $B_{S^k_j}$ is the greatest ball contained in $S^k_j$. We denote
by $X^k=\{x_i\}_{i\in\N}$ the corresponding set of the vertices,
 and introduce the space of continuous piecewise linear functions
on $\mathcal{T}^k$,
\[
W^k=\{w\in C(\R^N) :\text{$Dw(x)$ is constant in $S_j^k$}\}.
\]
Every element $w$ in $W^k$ can be expressed as
\begin{align*}
w(x)=\sum_{i\in \N}\beta_i(x)w(x_i),
\end{align*}
for basis functions (the so-called tent functions) $\beta_i\in W^k$
satisfying $\beta_i(x_j)=\delta_{ij}$ for $i,j\in\N$. It immediately
follows that $0\leq \beta_i(x)\leq 1$, $\sum_{i\in\N}\beta_i(x)=1$,
$\beta_i$ has compact support, and at any $x\in\R^N$ at most $N+1$
$\beta_i$'s are non-zero. The family $\{\b_j\}_j$ is a partition of
unity. On any simplex $S_i^k$, the $\beta_j$'s are called the
barycentric coordinates of $S_i^k$.

The fully discrete scheme can then be formulated as follows:
Find the function $u\in W^k$ that satisfies \eqref{HJBFh} at every
vertex $x_i\in X^k$, or equivalently,
\begin{equation}\label{fullydiscrete2}
    u(x_i)=\inf_{v\in V}\Big\{e^{-hc(x_i,v)} \Big[e^{-\l h}
    \sum_{j}M_{ij}(v)u(x_j)+(1-e^{-\l
    h})\sum_{j}P_{ij}(v)u(x_j)\Big]+hf(x_i,v)\Big\},
\end{equation}
for every $x_i\in X^K$. Here the matrices $M(v)$ and $P(v)$ are given
 by
\begin{align}
\label{Mv}
M(v)&=\sum_{m=1}^d \frac{1}{2d}( M_m^+(v)+M_m^-(v)) \\
\intertext{for $M^{\pm}_{m,ij}(v)=\beta_j\big(x_i+h
 b(x_i,v)\pm\sqrt{h}\s_m(x_i,v)\big),$ and}
 P_{ij}(v)&=\frac1\l\int_E \b_j(x_i+\eta(x_i,v,z))\n(dz).\label{Pv0}
\end{align}
Note that $M$ is a stochastic matrix, and for any $m$, only $N+1$
entries of any row of $M^{\pm}_m$ are non-zero. The matrix $P$ is
non-zero only if the vertex $x_j$ belong
to a simplex which has nonempty intersection with the set
$x_i+\eta(x_i,v,\mathrm{supp}(\nu))$ for all $v\in V$.

As a final step we also discretize $P_{ij}$ by (monotone)
quadrature
\begin{align}
    \label{Iappr}
    Q_{\dz}[\phi]:=\sum_{j\in\N}\phi(z_j)
    \omega_j \quad\text{where} \quad z_j\in E,\quad \omega_j\geq0,
\end{align}
satisfying the error bound
\begin{align}
    \label{int_err}
    E_{\dz}[\phi]:=\Big|\int_E\phi(z) dz -
     Q_{\dz}[\phi]\Big|\leq \bar K \|D\phi\|_{L^1}
     \dz\quad\text{for}\quad \phi\in W^{1,1}(\R^n).
\end{align}
Here $\dz$ is the discretization parameter.
All sensible tensor product quadratures satisfy this (first order)
error bound, and the bound is optimal if the integrand $\phi$ is not
more regular. The
monotonicity assumption $\omega_j\geq0$ is satisfied
for compound Gauss and Newton-Cotes types of quadratures in any
space dimension and when the order is less than 9 in the
Newton-Cotes case. We refer to \cite{JKLC} for examples and a wider
discussion of these issues. If we assuming that $\nu$ has a density
$m$ satisfying (B2), we get the following final scheme: Find $u\in
W^k$ such that
\begin{equation}\label{fullydiscrete}
    u(x_i)=\inf_{v\in V}\Big\{e^{-hc(x_i,v)} \Big[e^{-\l h}
    \sum_{j}M_{ij}(v)u(x_j)+(1-e^{-\l
    h})\sum_{j}\bar P_{ij}(v)u(x_j)\Big]+hf(x_i,v)\Big\},
\end{equation}
where $M_{ij}$ is defined in \eqref{Mv} and
\begin{equation*}
    \bar P_{ij}(v)=\frac1\l
    Q_{\dz}[\b_j(x_i+\eta(x_i,v,\cdot))m(\cdot)].
\end{equation*}
We have the following existence, uniqueness, and {\em partial} convergence
result for \eqref{fullydiscrete}.

\begin{thm}
\label{thm_dconv}
Assume (A1) -- (A4), (B2), \eqref{Iappr},
\eqref{int_err} hold, $h\l\leq \bar C_0$, and $h\in(0,1]$. Then
there exists a unique bounded  solution $u_{hk\dz}\in W^k$ to
\eqref{fullydiscrete}. Furthermore, if the solution $u_h$ of
\eqref{HJBFh} belongs to $C^{0,1}(\R^N)$, then
\begin{equation*}
    |u_h-u_{hk\dz}|_0\le \frac{|u_h|_1}{1-e^{-c_0h}}\Big(2k+C\dz\Big).
\end{equation*}
\end{thm}
\begin{proof}
Existence and uniqueness follow from a fixed point argument like in
the proof of Proposition \ref{propscheme1}. To prove the error bound,
note that since
$u_{hk\dz}(x)=\sum_{j\in\N}\b_j(x) u_{hk\dz}(x_j)$, $\b_j\geq0$, and
$\sum\b_j=1$,
\begin{align}
    |u_h(x)-u_{hk\dz}(x)|\le&\sum_{j\in\N}\Big(\b_j(x)|u_h(x)-u_{h}(x_j)|+\b_j(x)|u_h(x_j)-u_{hk\dz}(x_j)|\Big)
    \label{estfd1}\\
    \le &\, |u_h|_{1}k +\sum_{j\in\N}\b_j(x)|u_h(x_j)-u_{hk\dz}(x_j)|.\nonumber
\end{align}
The last term can be estimated by using \eqref{HJBFh} and
\eqref{fullydiscrete}. Easy computations show that
\begin{align*}
&|u_h(x_j)-u_{hk\dz}(x_j)|\\
&\le e^{-h c_0} e^{-\l h}|u_{hk\dz}-u_h|_0+e^{-h
  c_0}(1-e^{-\l h})\bigg[\sum_j\bar P_{ij}|u_h(x_j)-u_{hk\dz}(x_j)
  |\\
&+\underset{  A}{\underbrace{\Big|\sum_j\bar P_{ij}u_h(x_j) -\sum_j
P_{ij}u_h(x_j)\Big|}}+\underset{   B}{\underbrace{
 \Big|\sum_jP_{ij}u_h(x_j)-\frac1\l\int_E
u_h(x_i+\eta(x_i,z,v))\nu(dz)\Big|}} \bigg] .
\end{align*}
Note that $\sum_j \bar P_{ij}=\frac1\l Q_{\dz}[1]=1$ by
\eqref{int_err} since $\sum\b_j=1$. Furthermore,
\begin{align*}
A&=\Big|\frac{1}{\l}\int_EI_ku_h(x_i+\eta(x_i,v,z))m(z)dz-Q_{\dz}\Big[I_ku_h(x_i+\eta(x_i,v,\cdot))m(\cdot)\Big]\Big|\\
&\leq \frac{\dz}{\l} \big\|D\big[I_ku_h(x_i+\eta(x_i,v,z))m(z)\big]\big\|_{L^1}\leq\frac{\dz}{\l}
\Big(|u_h|_1L_3C_0+|u_h|_0C_1\Big)\int_Ee^{-\eps|z|}dz,\\
B&=\Big|\frac{1}{\l}\int_E(u_h-I_ku_h)(x_i+\eta(x_i,v,z))\nu(dz)\Big|\leq
 k  |u_h|_1\frac {\int_E\nu(dz)}{\l}=k  |u_h|_1,
\end{align*}
where $I_k\phi(x)=\sum_i\beta_i(x)\phi(x_i)$ is piecewise linear
interpolation of $\phi$ on $X^k$ and we have used (B2), (A4), and
\eqref{int_err}.
Combining these estimates and \eqref{estfd1}, using the properties
of $\beta_i(x)$, and remembering that $h\in(0,1]$ then gives the result.
\end{proof}

\begin{rem}
Since $h\in(0,1]$, $|u_h-u_{hk\dz}|_0\le
\frac{C|u_h|_1}{c_0}\frac{k+\dz}{h}$, which is consistent with the
estimates obtained in \cite{Ca:IPDE} in a different setting.
\end{rem}

\begin{rem}
If we use the cut-off procedure explained in Section
\ref{sec:BndDom} (a way of reducing to a bounded domain), then
equation \eqref{fullydiscrete} gives a finite system of equations.
In view of the integral term the new effective domain then becomes
$B_\mu+\mathrm{supp}(\nu)$ (see Section \ref{sec:BndDom}).
\end{rem}
\subsection{Schemes for unbounded measures $\nu$}
\label{sec:unbnd}
We consider general unbounded Levy measures $\nu$ under assumption
(B3).
There are two different cases: (i) $\alp\in[0,1)$ with HJB
equation \eqref{HJBF}, and (ii) $\alp\in[1,2)$ with HJB
  equation \eqref{HJB}. To derive our schemes,  we first reduce
to a bounded Levy measure $\nu_{r,R}$ as explained in Section
\ref{sec:prelim}. The result are the approximate HJB equations
\eqref{T1HJB} and \eqref{THJB}. These equations are then approximated
by a slightly modified version of the semi-Lagrangian scheme \eqref{HJBFh} (or
equivalently \eqref{DPPh}) defined in Section \ref{sec:timediscr}.

We propose the following semi-Lagrangian scheme in case (i)
\begin{align}\label{HJBFh22}
&v_h(x)=\inf_{v\in V }\Big\{hf(x,v)+
 e^{-h c(x,v)}\Big[\frac{e^{-h \l_{r,R} }}{4d}\sum_{m=1}^d \Big(v_h(x+h
   \bar b(x,v)+\sqrt{h}\bar\s_{+,m}(x,v))\\
&+v_h(x+h
   \bar b(x,v)+\sqrt{h}\bar\s_{-,m}(x,v))+v_h(x+h \bar
 b(x,v)-\sqrt{h}\bar\s_{+,m}(x,v))\nonumber\\
&+v_h(x+h
   \bar b(x,v)-\sqrt{h}\bar\s_{-,m}(x,v))\Big)+
\frac{1-e^{-h \l_{r,R}}}{\l_{r,R}}\int_E
v_h(x+\eta(x,v,z))\nu_{r,R}(dz)\Big]\Big\},\nonumber\\
\intertext{and in case (ii)}
\label{HJBFh33}
&w_h(x)=\inf_{v\in V }\Big\{hf(x,v)+
 e^{-h c(x,v)}\Big[\frac{e^{-h \l_{r,R} }}{4d}\sum_{m=1}^d (w_h(x+h
   \tilde b(x,v)+\sqrt{h}\bar\s_{+,m}(x,v))\\
&+w_h(x+h
   \tilde b(x,v)+\sqrt{h}\bar\s_{-,m}(x,v))+w_h(x+h \tilde
 b(x,v)-\sqrt{h}\bar\s_{+,m}(x,v))\nonumber\\
&+w_h(x+h \tilde b(x,v)-\sqrt{h}\bar\s_{-,m}(x,v))+ \frac{1-e^{-h
\l_{r,R}}}{\l_{r,R}}\int_E
w_h(x+\eta(x,v,z))\nu_{r,R}(dz)\Big]\Big\},\nonumber
\end{align}
where $\l_{r,R}:=\int_E\nu_{r,R}(dz)$, $\bar b, \tilde b, \nu_{r,R}$
are defined in Section \ref{sec:prelim}, and $\bar\s_{\pm,m}$ is
$m$-th column of
\begin{align}
\label{bspm}
\bar\s_{\pm}(x,v)=\s(x,v)\pm
\sqrt{\int_{0<|z|<r}\eta(x,v,z)\eta(x,v,z)^T\nu(dz)},
\end{align}
where the square root denotes the matrix square root.

\begin{rem}
The additional terms in \eqref{HJBFh22} and \eqref{HJBFh33} compared
with \eqref{HJBFh}, enable us to use $\s\pm\sqrt{\int\eta\eta^T\nu}$
instead of $\sqrt{\s\s^T+\int\eta\eta^T\nu}$ as diffusion matrix. The
consistency relation Proposition \ref{propscheme1} (v) still holds,
and if $\eta(x,v,z)=\eta_1(x,v)\eta_2(z)$, then the square root
in \eqref{bspm} equals $C\eta_1(x,v)$ where $C$ is the precomputable
constant matrix $\sqrt{\int_{0<|z|<r}\eta_2(z)\eta_2^T(z)\nu(dz)}$.
\end{rem}
\begin{rem}
 These schemes are similar to \eqref{HJBFh}, and can be
derived in a similar way. The conclusions of Propositions \ref{propscheme1} and
\ref{ContLipSch} (when $b$-terms are replaced by $\bar b$- or
$\tilde b$-terms as defined in Section \ref{sec:prelim}) still hold
for solutions of \eqref{HJBFh22} and \eqref{HJBFh33}. We refer to
\cite{JKLC} for the technical modifications needed to handle the integral term
in the diffusion coefficients.
\end{rem}

\begin{rem}
Previously bounded quantities may blow up as $r\ra0$. Indeed by
 (B3)  and (A2) we have for $r\in(0,1)$,
\begin{align}
\label{blowup}
&\l_{r,R} \leq \frac
{K}{\alp r^{\alp}} \quad \text{for } \alp\in(0,2),\qquad
\Big|\int_E\eta\nu_{r,R}(dz)\Big|\leq
 \left\{\begin{array}{ll}
 \frac{L_2K}{(\alp-1) r^{\alp-1}} &\text{for } \alp\in(1,2)\\
 L_2K & \text{for } \alp\in(0,1).
 \end{array}\right.
\end{align}
\end{rem}

\section{Convergence estimates for the discrete-time problem}
\label{sec:conv}

In this section we prove a priori error bounds for the convergence of
solutions $u_h$ of the semi-discrete scheme \eqref{HJBFh} to the
unique viscosity solution $u$ of \eqref{HJBF}. We consider two cases:
(i) The measure $\nu$ is bounded and (B1)$_0$ holds and (ii) the measure
$\nu$ is the truncation of an unbounded measure satisfying (B3).

In view of the equi-boundedness, equi-continuity  and consistency
results of Propositions \ref{propscheme1} and \ref{ContLipSch},
$u_h$ converge locally uniformly to $u$ by the Arzela-Ascoli
Theorem, stability and uniqueness results for viscosity solutions
(see e.g. \cite{BI07}), and the consistency result in Proposition
\ref{propscheme1}. It is also possible to obtain convergence without
equi-continuity (i.e. under weaker assumptions on the coefficients)
using so-called half relaxed limits \cite{BaSo}. Such results are
given in \cite{BLCN} for some non-local equations, but these results
does not cover the HJB equations we consider here.

Now we proceed to obtain a priori estimates for the convergence of
$u_h$ to $u$. To do this we will make use of an abstract result in
\cite{JKLC}, which we will describe below. Consider the equation
\begin{equation}\label{Continuous}
    F(x,u,Du, D^2 u ,u(\cdot))=0\qquad x\in\R^N
\end{equation}
where $u(\cdot)$ represents non-local (integral) terms.  Let $h>0$ be a
discretization parameter and consider an approximation scheme for
\eqref{Continuous} written in abstract form as
\begin{equation}\label{Scheme}
    S(h,x,u_h(x),[u_h]_x)=0\qquad x\in \R^N,
\end{equation}
where $[u_h]_x$ represents a function defined at $x$ via all the possible
value of $u_h$.  We need the following set of assumptions.
\begin{itemize}
\item[\textbf{(C1)}] \textit{(Monotonicity)} There exists $\ov c_0>0$
such that for any $h>0$, $x\in\R^N$, $\zeta\in \R$, $\t>0$ and
bounded functions $u$, $v$ such that $u\le v$ in $\R^N$, then
\[
 S(h,x,\zeta+\t,[u+\t]_x)\ge  S(h,x,\zeta,[v]_x)+\ov  c_0 \t
 \]
\item[\textbf{(C2)}] \textit{(Regularity)}
For any $h>0$ and any continuous, bounded function $\phi$, the function
\[x\mapsto S(h,x,\phi(x),[\phi]_x)\]
is bounded and continuous on $\R^N$ and the function
\[\zeta\mapsto S(h,x,\zeta,[\phi]_x)\]
is uniformly continuous for bounded $\zeta$, uniformly in $x$.
\medskip
\item[\textbf{(C3)}] {\em (Consistency)} There exists a function
  $E(\tilde K,h,\e)$ such
that for any sequence $\{\phi_\eps\}_\eps$ of smooth functions satisfying
$$  |D^{\beta} \phi_\e (x) | \leq \tilde K
\e^{1-|\beta|}  \quad \hbox{in  }\R^N ,\quad\text{ for
any  $\beta\in\N^{N}$},$$
where $|\beta|=\sum_{i=1}^N \beta_i$, the following inequality holds:
$$ |S(h,x,\phi_\e (x),[\phi_\e
]_{x})-F(x,\phi,D\phi_\e,D^2\phi_\e)|\leq E (\tilde K,h,\e)\quad
\hbox{in  }\R^N.$$

\item[\textbf{(C4)}] \textit{(Convexity)} Let $(\rho_\eps)_{\eps>0}$
be a family of mollifiers (smooth, positive functions with mass 1 and
support in $\{|x|<\eps\}$).  For any Lipschitz-continuous function
$\phi$, there exists a constant $C$ such that for any $x$ and $h$
\[\int_{\R^N} S(h,x,\phi(x-e),[\phi(\cdot-e)]_x)\rho(e)de\ge
S(h,x,(\phi*\rho_\eps)(x),[(\phi*\rho_\eps)]_x)-C\eps
\]
\item[\textbf{(C5)}] \textit{(Commutation with translations)} For any
$h>0$ small enough, $0<\eps<1$, $y\in \R^N$, $\zeta\in\R$, continuous
bounded function $\phi$ and $| e|\le y$, we have
\[
S(h,y,\zeta,[\phi]_{y-e})=S(h,y,\zeta,[\phi(\cdot-e)]_y).
\]
\item[\textbf{(D)}] For $h$ small enough and $\eps\in[0,1)$, there is a
  unique solution $u_h^\eps$ of the scheme
\begin{equation}\label{HJBFhe}
\max_{|e|\le \eps} S(h,x+e,u^\eps_h(x),[u^\eps_h]_{x})=0\quad\text{in } \R^N,
\end{equation}
where $u_h:=u_h^0$ also solve \eqref{Scheme}, and a constant $C$
independent of $h,\eps$ such that
\[
|u_h^\eps|_1\le C\,\qquad\text{and} \qquad
|u_h^0-u_h^\eps|_0\le C\eps.
\]
\end{itemize}

We remark that we are using a more general consistency relation here than
in \cite{JKLC}, and that this extra generality will be needed when we
consider unbounded measures $\nu$. The next result is a restatement of
Theorem 3.4 in \cite{JKLC} in view of the new consistency relation (C3).

\begin{thm}
\label{abstract}
Assume (A1) -- (A3), (B1)$_2$, (C1) -- (C5), and (D) hold, and let $u$ and
$u_h$ be solutions of respectively \eqref{Continuous} and
\eqref{Scheme} satisfying $\tilde K:=|u|_1\vee |u_h|_1<\infty$.
Then there exists a constant $C$ depending only on $L_1, L_2, c_0$
from (A2) and (A3) such that
$$| u-u_h |\leq C\min_{\e>0} \left(\e + E_1 (\tilde K,h,\e)\right)
\quad\text{in}\quad\R^N.
$$
\end{thm}

\begin{rem}
To prove (D) for the scheme \eqref{HJBFh}, we will need to assume also
(B1)$_0$, $h\lambda\leq \bar C_0$, $h\leq 1$, and $c_0\geq c_1$ for
both $c_1$ defined in Theorem \ref{WP} and Proposition
\ref{ContLipSch}. Under these assumptions we also have $\tilde
K:=|u|_1\vee |u_h|_1<\infty$.
\end{rem}

We will apply this abstract result to derive error bounds for the scheme
\eqref{HJBFh}. We rewrite the scheme in the form \eqref{Scheme} with
\begin{align}
&S(h,x,r, [\psi]_x)=\nonumber\\
&\sup_{v\in V}\Big\{ \frac{-e^{-\l
h}}{2dh}\sum_{m=1}^d\Big[[\psi]_x(h
b(x,v)+\sqrt{h}\s_m(x,v))-2r+[\psi]_x(hb(x,v)
-\sqrt{h}\s_m(x,v))\Big]\label{S}\\
&-\frac{1-e^{-\l h}}{\l
h}\int_E[\psi(x +\eta(x,v,z))-r]+\frac{e^{hc(x,v)}-1}{h} r-e^{hc(x,v)}
f(x,v)\Big\},\nonumber
\end{align}
and $[\psi]_x(z)=\psi(x+z)$.
\begin{lem}\label{propscheme2}
Assume (A1) -- (A3) and (B1)$_0$, $h\l\leq \bar C_0$, $h\le1$. Then the scheme
\eqref{HJBFh} (and equivalently \eqref{S}) satisfies assumptions
(C1)--(C5) with
$$E(\tilde
K,h,\eps)=C_1h(\tilde
K\eps^{-1}+\tilde
K\eps^{-2}+\tilde
K\eps^{-3})
+C_2h\lambda((1+|\int_E\eta\nu|)\tilde
K+\tilde
K\eps^{-1}),$$
where the constant $C_1$ and $C_2$ only depend on
$\sup_v|\sigma|_0,\sup_v|b|_0$.

If in addition $c_0\geq c_1$ for both $c_1$'s in Theorem \ref{WP}
and Proposition \ref{ContLipSch}, then assumption (D) also holds
with constants $C$ only depending on the data and $\bar C_0$.
\end{lem}
\begin{proof}
It is straightforward to verify (C1) with $\ov c_0=c_0$ where $c_0$
is defined in (A3). (C2) follows from the assumptions on the
coefficients, while (C3) follows from Proposition \ref{propscheme1}
(v). By a straight forward computation, it follows that (C4) holds
with $C=0$. We refer to \cite{JKLC} for similar computations.
Finally, (C5) holds since $[\phi]_{x-e} =[\phi(\cdot-e)]_x$  for any
continuous function $\phi$.

To prove (D), observe that \eqref{HJBFhe} can be rewritten
in the form \eqref{S} (by defining a new control $\bar v=(v,e)$). The
coefficients of this new equation still satisfy (A1) - (A3). Therefore
(D) follows after an application of Propositions \ref{propscheme1}
and \ref{ContLipSch}.
\end{proof}

Now we are in a position to state the error bounds. First we
consider the bounded case, i.e (B1)$_0$ holds. In this case the
equation is \eqref{HJBF} which is approximated by the scheme
\eqref{HJBFh}. Note that the integral operator has not yet been
discretized.

\begin{thm}[Bounded measure]
\label{main_bnd} Assume (A1) -- (A3), (B1)$_0$, $h\le1$, and
$c_0\geq c_1$ for both $c_1$'s in Theorem \ref{WP} and Proposition
\ref{ContLipSch}. Let $u$ be the solution of \eqref{HJBF} and $u_h$
be the solution of \eqref{HJBFh} (or equivalently \eqref{S}).

\smallskip
(a) (General IPDEs) Then $|u-u_h|_0\leq Ch^{1/4}.$
\smallskip

(b) (1st order IPDEs) If in addition $\sigma\equiv0$, then
$|u-u_h|_0\leq Ch^{1/2}.$
\smallskip

\noindent In both cases the constant $C$ depends only on the
coefficients, $c_1$, and $\l$.
\end{thm}

\begin{proof}
Part (a) is an easy consequence of Theorem \ref{abstract} and
Lemma \ref{propscheme2}. Part (b) follows in a
similar manner after noting that the consistency relation
corresponding to Proposition \ref{propscheme1} {\it (v)} now becomes
\begin{align*}
&|L^v\psi(x) - L_h^v\psi(x)|\leq
C_1h|D^2\phi|_0+C_2h\lambda((1+|\int_E\eta\nu|)|D\phi|_0+|D^2\phi|_0).
\end{align*}
\end{proof}
\begin{rem}The convergence rate obtained in Theorem
\ref{main_bnd} is the same as in the pure PDE case, see
\cite{J,BJ:Err1}.  In the first order case and when the Levy measure
is bounded, convergence rate 1/2 for semi-Lagrangian schemes like
\eqref{HJBFh} have previously be obtained in \cite{Ca:IPDE}. However
the integral term in \cite{Ca:IPDE} has a different form compared to
the one we consider here.
\end{rem}

\begin{rem}
The scheme \eqref{HJBFh} uses a first order accurate approximation
of 2nd derivatives as can be seen from the consistency relation
Proposition \ref{propscheme1} {\it (v)}. This leads to lower rates
of convergence than for some monotone FDMs that use 2nd order
accurate approximations of 2nd derivatives, see
\cite{Kry05,JKLC,BJK2}. There the rate is 1/2, while for a more
general class of FDMs the rate is at least 1/5, see \cite{BJ:Err2}
for the pure PDE case.
\end{rem}

We proceed to the case of general unbounded Levy measure $\nu$ under
assumption (B3). There are two different cases: (i) $\alp\in[0,1)$
with HJB equation \eqref{HJBF}, and (ii) $\alp\in[1,2)$ with HJB
  equation \eqref{HJB}.
In case (i), $\alp\in(0,1)$, as a consequence of Lemma
\ref{propscheme2}, Theorems \ref{abstract} and \ref{TruncThm} we have
the following convergence result.

\begin{thm}[Unbounded measure I]
\label{main_unbnd1} Assume (A1) -- (A3), (B3) with $\alp\in(0,1)$,
$h\le1$, and $c_0\geq c_1$ for both $c_1$'s in Theorem \ref{WP} and
Proposition \ref{ContLipSch}. Let $u$ be the solution of
\eqref{HJBF} and $u_h$ be the solution of \eqref{HJBFh22}.

Then the
best rate is obtained choosing $r=h^{\frac3{6+\alp}}$, and in this case
$$|u-u_h|_0\leq C(h^{1/4}+h^{\frac{3-\alp}{6+\alp}})\leq Ch^{1/4},$$
where the constant $C$ depend only on the coefficients, $c_1$, and
quantities from (B3)/(B1)$_2$.
\end{thm}

\begin{proof}
When $\alp\in(0,1)$, Lemma \ref{propscheme2} still holds for the
scheme \eqref{HJBFh22}, and in view of \eqref{blowup}
we have the following form of $E$,
$$E(\tilde
K,h,r,\eps)=C_1\tilde K h(\eps^{-1}+\eps^{-2}+\eps^{-3})+C_2\tilde K
h r^{-\alp}(1+\eps^{-1}),$$ where the constant $C_1$ and $C_2$ are
independent of $h,r,\eps$. Let $u_r$ denote the solution of
\eqref{T1HJB}. Theorem \ref{abstract} and (``term-wise'')
minimization in $\eps$, lead to the bound
$$|u_r-u_h|_0\leq C(h^{1/4}+r^{-\alp}h+r^{-\frac\alp2}h^{\frac12}),$$
where the constant $C$ depend only on the coefficients, $c_1$, and
quantities from (B3)/(B1)$_2$. In view of Theorem \ref{TruncThm} and the optimal
choice of $r$, $r=h^{\frac3{6+\alp}}$, the result follows.
\end{proof}

\begin{rem} Since $\frac14<\frac{3-\alp}{6+\alp}<\frac12$ for
  $\alp\in(0,1)$, there is no reduction of rate due to truncation of
  the measure $\nu$.
\end{rem}

The case $\alp\in(1,2)$ is more difficult, since now also
$\int_E\eta\nu_{r,R}$ and hence $\tilde b$ in \eqref{HJBFh33} blows up
as $r\ra0$. As a consequence Theorem \ref{abstract} can no longer be
used directly. The convergence result is the following:

\begin{thm}[Unbounded measure II]
\label{main_unbnd2} Assume (A1) -- (A3), (B3) with $\alp\in(1,2)$,
$h\le1$, and $c_0\geq c_1$ for both $c_1$'s in Theorem \ref{WP} and
Proposition \ref{ContLipSch}. Let $u$ be the solution of \eqref{HJB}
and $u_h$ be the solution of \eqref{HJBFh33}.

Then the best rate is obtained choosing $r=h^{\frac3{3+5\alp}}$, and
in this case
$$|u-u_h|_0\leq Ch^{\frac{3-\alp}{3+5\alp}},$$
where the constant $C$ depend only on the coefficients, $c_1$, and
quantities from (B3)/(B1)$_2$.
\end{thm}
\begin{rem}
This result is consistent with Theorem \ref{main_unbnd1} since
  $\underset{{\alp\ra1^+}}{\lim}\frac{3-\alp}{3+5\alp}=\frac14$. For $\alp\in(1,2)$
  the rate degrades as $\alp$ increases, and
  $\underset{\alp\ra2^-}\lim\frac{3-\alp}{3+5\alp}=\frac1{13}$.
\end{rem}
\begin{proof}[Outline of proof]
 From Proposition \ref{ContLipSch} we have
uniform in $r$ Lipschitz continuity, but the continuous dependence
estimates will be proportional to $\int_E|z|\nu_{r,R}$ through the
$\tilde b$-term. Because of this we must redo the arguments leading to
Theorem \ref{abstract}, and the result will be an estimate of the form
$$
| u_r-u_h |\leq C\min_{\e>0} \left(\e\int_E|z|\nu_{r,R} + E_1
(\tilde K,h,r,\e)\right) \quad\text{in}\quad\R^N,
$$
when $u_r$ solve \eqref{THJB}. We omit the details, since the
argument is same as the one used to prove Theorem \ref{abstract} in
\cite{JKLC}. Because of the blow up in $\int_E\eta\nu_{r,R}$, we
also need a much more precise consistency relation than given in
Proposition \ref{propscheme1} (v) tracking all $\tilde b$
dependence. From the proof of Proposition \ref{propscheme1}, it is
easy to see that it will have the following form
\begin{align*}
E(\tilde K,h,\eps)= &\ C_1\tilde
K\Big[h|\tilde b|_0^2\eps^{-1}+(h^2|\tilde b|_0^3+h
|\tilde b|_0)\eps^{-2}+(h+h^3|\tilde b|_0^4+h^2|\tilde b|_0^2)\eps^{-3}\Big]\\
&+C_2\tilde
K\l_r h\Big[|\tilde b|_0+\eps^{-1}\Big]+C_3\tilde
K\l_r h\Big[\int_E\eta\nu_{r,R}\Big],
\end{align*}
where the constants $C_1, C_2, C_3$ only depend on $\sup_v|\bar
\sigma|_0$. In view of \eqref{blowup},
$\l_r=Cr^{-\alp}$ and $|\tilde b|_0+\int_E|z|\nu_{r,R}\leq
C(1+r^{\alp-1})$. The rest of the proof is a long computation
consisting of choosing optimal $\eps$ and $r$ as in the proof of
Theorem \ref{main_unbnd1}. We omit the details only remarking that
the worst term in $E$ turns out to be the $h\eps^{-3}$-term.
\end{proof}

\appendix

\section{Lipschitz regularity and continuous dependence}
\label{sec:CDP}

In this section we obtain a combined Lipschitz regularity and
continuous dependence result for solutions to the equation
\begin{align}
u_h(x)&=\inf_{v\in V }\Big\{hf(x,v)+ e^{-h c(x,v)}\Big[\frac{e^{-h \l }}{2d}\sum_{m=1}^d \Big(u_h(x+h\ov b(x,v)+\sqrt{h}\s_m(x,v))+\label{EQ}
\\
&+u_h(x+h\ov b(x,v)-\sqrt{h}\s_m(x,v)\Big)+
\frac{1-e^{-h \l}}{\l}\int_E u_h(x+\eta(x,v,z))\nu(dz)\Big]\Big\},\nonumber
\end{align}
where $\lambda=\nu(E)$ and
$$\bar b=b+\frac{1-e^{-h\lambda}}{\lambda e^{-\lambda h}}\int_E\eta \nu.$$

\begin{thm}
\label{LipCont} Let $u_h$ and $\tilde u_h$ be solutions of
\eqref{EQ} corresponding to the data $\sigma, \bar b, c, f, \eta,
\nu$ and $\bs, \bar\bb, \bc, \bff, \bbeta, \nu$ respectively, and
assume both sets of coefficients satisfy (A1) -- (A3) and (B1)$_0$,
and that $h\lambda\leq C$ and $h\leq 1$. Then there exist constants
$c_1, L, K\geq0$ (only depending on the data and $C$) such that if
$c_0\geq c_1$, then for all $h>0,x,y\in\R^N$,
\begin{align*}
|u_h(x)-\bu_h(y)|\leq & \ L|x-y|
+K\sup_{v\in\V}\big[|f-\bff|_0+|c-\tilde
c|_0\\
&+|b-\bb|_0+|\sigma-\bs|_0+
|\int|\eta(\cdot,z)-\bbeta(\cdot,z)|^2\nu(dz)|_0^{1/2}\big].
\end{align*}
\end{thm}

\begin{rem}
The precise dependence of the constants $c_1,L,K$ is given in the
proof below.
\end{rem}

\begin{rem}
This result extends the corresponding results of \cite{BJ:Err1} to
non-local HJB equations, with general (singular) Levy
measures. Moreover, the proof below simplifies the corresponding
proofs of Barles and Jakobsen \cite{BJ:Err1}
because we do not use the somewhat unnatural ``doubling schemes'' as
in \cite{BJ:Err1}. Instead we give a more direct proof.
\end{rem}

\begin{proof}
We will use doubling of variables techniques similar to those
used to prove corresponding results for equation \eqref{HJBF}.
We define
\begin{align}
    \phi(x,y)&=m+\alp M+L(\alp^{-1}+\alp|x-y|^2)+\eps(|x|^2+|y|^2),\\
    \psi(x,y)&=\sup_{x,y}[u_h(x)-\bu_h(y) - \phi(x,y)] = \psi(\bx,\by),
\end{align}
where $\alp,\eps>0$, $m,M,L\geq 0$, and $(\bx,\by)$ is the point where
the supremum is attained. We will prove that $\psi(\bar
x, \bar y)\leq o(1)$ as $\eps\rightarrow0$ for a suitable constant
$L$, and for
\begin{align}
    m&=\Big(\frac{e^{c_0h}-1}{h}\Big)^{-1}\sup_{v\in\V}e^{|c|_0\vee |\tilde c|_0
        h}\Big[|f-\bff|_0+(|u_h|_0\wedge|\tilde
    u_h|_0 + h |f|_0\wedge |\tilde f|_0)|c-\tilde c|_0\Big],\label{m}\\
    M&= K\Big(\frac{e^{c_0h}-1}{h}\Big)^{-1}\sup_{v\in\V}
    \big[|\sigma-\bs|_0^2+|b-\bb|_0^2+|\int|\eta(\cdot,z)-\bbeta(\cdot,z)|^2\nu(dz)|_0\big]. \label{M}
\end{align}
This implies Theorem \ref{LipCont} (after sending $\eps\ra 0$) because for any
$\alp\in\R$ and $x,y\in\R^N$,
$$u_h(x)-\bu_h(y)-m-3(LM)^{1/2}-L|x-y|\leq \psi(x,y)\leq \psi(\bar x,
\bar y)\leq o(1) \quad\text{as } \eps\ra0.$$

For simplicity we will not be explicit about the form of the
$\eps$-terms appearing in the computations below. The role of these terms is
only to guaranty that the maximum is attained at a point $(\bx,\by)$, and
they vanish in the final estimate when $\eps\ra0$. We refer to the
proof of Theorem 3.4 in \cite{BJ:Err1} for details concerning the $\eps$-terms.

We proceed by contradiction assuming that $\psi(\bx,\by)> o(1)$ as
$\eps\ra0$. Note that by the definition of $\psi$ this implies that
$u_h(\bx)-u_h(\by)>0$. Furthermore, observe that since $(\bx,\by)$ is
a maximum point,
\begin{align*}
    \psi(\bx+b+a,\by+\bar b+\bar a)+\psi(\bx+b-a,\by+\bar b-\bar a)&\leq 2
    \psi(\bx,\by),\\
    \psi(\bx+\zeta,\by+\bar\zeta)&\leq \psi(\bx,\by),
\end{align*}
for every $a,b,\zeta,\bar a, \bar b,\bar\zeta\in\R^N$, and hence by
the definition of $\psi$,
\begin{align*}
    &I_1:=[u_h(\bx+b+a)-2v_h(\bx)+u_h(\bx+b-a)]\\
    &\qquad\quad-[\bu_h(\by+\bar b+\bar
        a)-2\bu_h(\by)+\bu_h(\by+\bar b-\bar a)]\\
    &\quad \leq \phi(\bx+b+a,\by+\bar b+\bar a)-2
        \phi(\bx,\by)+\phi(\bx+b-a,\by+\bar b-\bar a),\\
    &I_2:=[u_h(\bx+\zeta)-u_h(\bx)]-[\bu_h(\by+\bar\zeta)-\bu_h(\by)]\leq
        \phi(\bx+\zeta,\by+\bar\zeta)-
    \phi(\bx,\by).
\end{align*}
Finally, by the definition of $\phi$ we are lead to
\begin{align}
\label{I1}
& I_1 \leq 2L\alp|a-\bar a|^2 + 2L\alp|b-\bar b|^2+4L\alp (\bx-\by,b-\bar b)+o(1) && \text{as }\eps\ra0,\\
\label{I2}
&I_2\leq 2L\alp(\bx-\by,\zeta-\bar\zeta) + L\alp |\zeta-\bar\zeta|^2 + o(1) &&
\text{as }\eps\ra0.
\end{align}
These two inequalities are crucial for the rest of the proof.

Now we divide \eqref{EQ} by $h e^{-h c}$ and rewrite it as
\begin{align*}
    &0=\sup_{v\in V }\Big\{\frac{e^{hc}-1}h u_h(x)- e^{h c}f(x,v)\\
    &-\frac{e^{-h \l }}{2d h}\sum_{m=1}^d \Big(u_h(x+h\ov
    b(x,v)+\sqrt{h}\s_m(x,v))-2u_h(x)+u_h(x+h\ov
    b(x,v)-\sqrt{h}\s_m(x,v)\Big)\\
    &-\frac{1-e^{-h \l}}{\l h}\int_E
        \Big(u_h(x+\eta(x,v,z))-u_h(x)\Big)\nu(dz)\Big\}.\nonumber
\end{align*}
Upon subtracting the equations (in this new form) for $\bu_h$ and
 $u_h$, we get
\begin{align*}
    &0\leq \sup_{v\in V }\Big\{\frac{e^{\tilde
        ch}-1}h\bu_h(\by)-\frac{e^{ch}-1}h u_h(\bx)\\
    &\qquad +[e^{ch}f(\bx,v)-e^{\tilde ch}\bff(\by,v)]+ \frac{e^{-h \l
        }}{2dh}\sum_{m=1}^d I_{1,m}+\frac{1-e^{-h \l}}{\l h}\int_E
        I_2\nu(dz)\Big\},
\end{align*}
where $I_{1,m}$ corresponds to $I_1$ above with the choice
$a=\sqrt{h}\s_m(\bx,v)$, $b=h\ov b(\bx,v)$, $\bar a=\sqrt{h}\bs_m(\by,v)$, and
$\bar b=h\bar\bb(\by,v)$, and for $I_2$ we have taken
$\zeta=\eta(\bx,v,z)$ and $\bar\zeta=\eta(\by,v,z)$. By \eqref{I1} and
\eqref{I2} we then have
\begin{align*}
    &0 \leq \sup_{v\in V }\Big\{\frac{e^{\tilde
        ch}-1}h\bu_h(\by)-\frac{e^{ch}-1}h u_h(\bx)\\
    &\quad +[e^{ch}f(\bx,v)-e^{\tilde ch}\bff(\by,v)]+\frac{e^{-h
        \l }}{2h}L\alp \Big[\frac1{d}\sum_{m=1}^d
    2  h|\s_m(\bx,v)-\bs_m(\by,v)|^2\\
    &\quad +2 h^2|\bar b(\bx,v)-\bar
    \bb(\by,v)|^2+4 h (\bx-\by, \bar b(\bx,v)-\bar \bb(\by,v))\Big]\\
    &\quad +e^{-h c}\frac{1-e^{-h \l}}{\l h}L\alp\int_E
    \Big[|\eta(\bx,v,z)-\bbeta(\by,v,z)|^2+
    2\big(\bx-\by,\eta(\bx,v,z)-\bbeta(\by,v,z)\big)\Big]\nu(dz)\Big\}\\
    &\quad+ o(1) \quad\text{as}\ \eps\ra 0.
\end{align*}
Since $\bar b=b+\frac{1-e^{-h\lambda}}{\lambda e^{-\lambda
h}}\int_E\eta \nu$ and $\bar\bb$ is defined similarly we see that
the $(\bx-\by,\eta-\bbeta)$-terms cancel in the above inequality.
Since  $\nu(E)=\lambda$, Jensen's inequality implies that
$$
    \lambda^2\Big|\int_E(\eta(\bx,v,z) -\bbeta(\by,v,z))
    \frac{\nu(dz)}{\lambda}\Big|^2\leq \lambda^2
    \int_E|\eta(\bx,v,z) -\bbeta(\by,v,z) |^2\frac{\nu(dz)}{\lambda}.
$$
Also note that since $u_h(\bx)-\tilde u_h(\by)>0$ and $c\geq c_0>0$,
$$
    \frac{e^{\tilde ch}-1}h\bu_h(\by)-\frac{e^{ch}-1}h u_h(\bx)\leq
    -\frac{e^{c_0h}-1}h [u_h(\bx)-\tilde u_h(\by)] + \frac{|e^{\tilde
        ch}-e^{ch}|}h \tilde |u_h(\bx)|\wedge |u_h(\by)|.
$$
Therefore after cancellations, Jensen's inequality, and the inequality $\frac{1-e^{-h\lambda}}{\lambda e^{-\lambda
    h}}\leq he^{h\lambda}$, we get
\begin{align*}
    &\frac{e^{c_0h}-1}{h}[u_h(\bx)-\bu_h(\by)]\\
    &\leq \sup_{v\in V
    }\Big\{e^{|c|_0\vee |\tilde c|_0
        h}\Big[|f(\bx,v)-\bff(\by,v)|+(|u_h|_0\wedge|\tilde
    u_h|_0 + h |f|_0\wedge |\tilde f|_0)|c(\bx,v)-\tilde c(\by,v)|\Big]\\
    & \quad + \frac{e^{-h
        \l }}{2h}L\alp\Big[\frac1{d}\sum_{m=1}^d
    2 h|\s_m(\bx,v)-\bs_m(\by,v)|^2\\
    &\quad +4 h^2|b(\bx,v)- \bb(\by,v)|^2+4 h (\bx-\by,  b(\bx,v)-
    \bb(\by,v))\Big]\\
    &\quad+\frac{1-e^{-h \l}}{\l h}L\alp\int_E \Big[4 h^2 \lambda he^{\lambda
        h}\int_E|\eta(\bx,v,z) -\bbeta(\by,v,z) |^2\nu(dz)\\
    &\quad + |\eta(\bx,v,z)-\bbeta(\by,v,z)|^2\Big]\nu(dz)\Big\}+ o(1) \quad\text{as}\ \eps\ra 0.
\end{align*}
Now since $\lambda h\leq C$ and $h<1$, it follows from simple
computations that
\begin{align*}
    &\frac{e^{c_0h}-1}{h}[\psi(\bar x,\bar y)+\phi(\bar x,\bar
    y)]=\frac{e^{c_0h}-1}{h}[u_h(\bx)-\bu_h(\by)]\\
    &\leq \sup_{v\in V }\Big\{\underset{\frac{e^{c_0h}-1}{h}m}{\underbrace{e^{|c|_0\vee |\tilde c|_0
        h}\Big[|f-\bff|_0+(|u_h|_0\wedge|\tilde
    u_h|_0 + h |f|_0\wedge |\tilde f|_0)|c-\tilde c|_0\Big]}}\\
    &\quad +L\alp \frac{e^{-h
        \l }}{2}\Big[\frac1{d}\sum_{m=1}^d
    4  |\s_m-\bs_m|_0^2 +16 h|b- \bb|_0^2\Big]\\
    &\quad+ L\alp 5Ce^C | \int_E|\eta(\cdot,v,z) -\bbeta(\cdot,v,z)
    |^2\nu(dz)|_0\Big\}\\
    &\quad+ \sup_{v\in V }\Big\{\underset{\overline
    L}{\underbrace{e^{|c|_0\vee |\tilde c|_0
        h}\Big[L_f+(|u_h|_0\wedge|\tilde
    u_h|_0 + h |f|_0\wedge |\tilde f|_0)L_c\Big]}}(\alp^{-1}+\alp|\bx-\by|)\\
    &\quad+ L\alp |\bx-\by|^2\underset{L_0}{\underbrace{\Big(\frac{e^{-h
        \l }}{2}\Big[\frac1{d}\sum_{m=1}^d
    4 L_{\sigma}^2 +16 hL_b+8\Big]+5Ce^C  L_{\eta}^2\int_E
    |z|^2\nu(dz)\Big)}}\Big\}\\
    &\quad+ o(1) \quad\text{as}\ \eps\ra 0.
\end{align*}
Let $\overline L, L_0$ be defined as in the inequality above. If $c_0$
is so big that
$$\frac{e^{c_0h}-1}{h}- L_0>0,$$
and we choose $m,M$ as in \eqref{m} and \eqref{M} for $K$ big enough, and
$$L=\frac{\overline
  L}{\frac{e^{c_0h}-1}{h}- L_0},$$
then $\psi(\bar x,\bar y)\leq o(1)$ as $\eps\ra 0$ and the proof is complete.
\end{proof}

\end{document}